\documentclass{svjour3}
\usepackage{amscd,amsfonts,amssymb,amsmath,latexsym,array,xcolor,graphicx}
\usepackage{times}
\newcommand\F{\mbox{I\kern-2pt F}}

\newcommand\cE{{\cal E}}

\newcommand\cF{{\cal F}}

\newcommand\cL{{\cal L}}

\newcommand\cT{{\cal T}}

\newcommand{\RNumb}[1]{\uppercase\expandafter{\romannumeral #1\relax}}

\def\E{{\bf E}}
\def\P{{\bf P}}

\def\bbr{{\mathbb R}}

\newcommand\beq{\begin{equation}}
\newcommand\eeq{\end{equation}}
\newcommand\bea{\begin{eqnarray}}
\newcommand\eea{\end{eqnarray}}
\newcommand\bean{\begin{eqnarray*}}
\newcommand\eean{\end{eqnarray*}}

\begin{document}

\title{On Integral Equations in the Ruin Theory with Risky Investments}

\author{Viktor Antipov \and Yuri Kabanov}

\institute{Viktor Antipov \at Lomonosov Moscow State University and ``Vega'' Institute, Moscow, Russia\\
\email{stayaptichek@gmail.com}
\and Yuri Kabanov (corresponding author) \at Lomonosov Moscow State University, Moscow, Russia, and Laboratoire de Math\'ematiques, Universit\'e de Franche-Comt\'e, Besan\c{c}on, France\\
\email{ykabanov@unif-fcomte.fr}}

\date{Received: date / Accepted: date}

\titlerunning{On Integral Equations in the Ruin Theory with Risky Investments}
\authorrunning{V. Antipov, Yu. Kabanov}

\maketitle

\begin{abstract}
We consider the ruin problem for an insurance company investing in a risky asset with the price given by a geometric Brownian motion and with the business process evolving as a L\'evy process with one-sided jumps and having trajectories of bounded variation. We consider the cases of non-life insurance and annuity payments and find, under rather mild assumptions on the L\'evy measure, the exact asymptotics of the ruin probabilities when the initial capital tends to infinity. Our note is complementary to the recent papers by Promyslov and Antipov--Kabanov, where the business processes were modeled by compound Poisson processes with drift. The arguments use the approach based on a verification theorem.

\end{abstract}

\section{Introduction}

In the modern economy, insurance companies may, and are even legally obliged to, invest their capital reserves in financial markets. This situation makes obsolete the classical collective risk theory, where it is assumed that the capital reserve is kept in a bank account with zero interest rate. Intensive studies that started in the nineties resulted in the creation of a collective risk theory with risky investments. There is a great difference between the conclusions of the ``old'' and the ``new'' theories. It turns out that financial risk strongly influences the ruin probabilities: as the initial capital grows, they decay much more slowly than predicted by the classical models. To explain the difference, let us consider as an example an insurance company investing the capital reserve in a composite index evolving as a geometric Brownian motion. In the case of exponentially distributed claims, the ruin probabilities decay, at best, as a power function of the initial capital, see, e.g., the papers \cite{FrKP}, \cite{KP2016}, \cite{ABL}, \cite{ACT}, \cite{KPukh} {\it et al.} where the arguments were based on the asymptotic theory of linear differential equations. Moreover, ruin is imminent if the volatility is high, and this qualitative conclusion can be quantified. This result is in a striking contrast with the Cram\'er--Lundberg theory predicting, optimistically, that in the case of positive safety loading, the ruin probabilities decrease exponentially fast.

The results obtained in the modern ruin theory cover a plethora of models and were proven using various mathematical tools. The discussion of the current state of the art is far beyond the scope of our short note. We refer the interested reader to detailed discussions in \cite{KP2020} and \cite{AK2026}. We want to emphasize that the mathematical tools used in the majority of studies are rather involved: implicit renewal theory, \cite{Go91}, \cite{Paul-93}, \cite{Paul-G}, theory of viscosity solutions, \cite{Bel-Kab}, \cite{AK2025}, Tauberian theorems, \cite{Grandits2004}, \cite{A2025}. This makes it difficult to present the derivation of these results, which look simple, to students of actuarial schools. Fortunately, the accessibility of the arguments has substantially improved since the recent publications \cite{Prom} and \cite{AK2026} inspired by the work of Grandits \cite{Grandits2004}. In these papers, a fairly simple approach was developed, in the spirit of verification theorems of the optimal control. Asymptotic results in the ruin problems for models of annuity payments as well as of non-life insurance can be explained using only background from a standard syllabus. The approach of the mentioned papers covers the case of the Cram\'er--Lundberg-type models where the business process is given by a compound Poisson process with drift and one-sided jumps, while the price of a risky asset follows a geometric Brownian motion. It also allows one to compute numerically the constants in the asymptotic expansions.

The modern literature considers more sophisticated frameworks based on L\'evy processes. For this reason, it is natural to ask whether the approach based on the study of integral equations describing the evolution of the density of the ruin probability can be applied to more general models. In the present note we show that it can be extended, under rather weak additional assumptions, to the case where the business activity is modeled by L\'evy processes with one-sided jumps and trajectories of bounded variation. This is the main contribution of the present note.

\section{Model}

We are given a filtered probability space $(\Omega,\cF,{\bf F}=(\cF_t)_{t\ge 0},\P)$ with two independent processes, a Wiener process $W=(W_t)$ and a L\'evy process $P=(P_t)$ of {\bf bounded variation} with $P_t=ct+x*\mu_t$, that is, in traditional differential notation,
\begin{equation}
dP_t=c dt+\int x\mu(dt,dx),
\end{equation}
where the jump measure $\mu=\mu(dt,dx)$ is the Poisson random measure with compensator $q(dt,dx)=dt\Pi_P(dx)$.

Put $R_t:=at+\sigma W_t$. The stochastic exponential
$S_t:=\cE_t(R):=e^{\sigma W_t +(a-\sigma^ 2/2)t}$
is the geometric Brownian motion describing the price evolution of the risky asset, while $P$ models the business process.

In the sequel we use the notation $h(x):=xI_{\{|x|\le 1\}}$ and $\bar h(x):=xI_{\{|x|> 1\}}$; the symbol of a function as the argument of a measure means the integral of this function. Recall that we always have $\Pi_P(|x|^2\wedge 1)<\infty$ and $P$ is of bounded variation if and only if $\Pi_P(|h|)<\infty$.

In the classical models, the process $P$ is a compound Poisson process with drift. It is usually represented in the following easily interpreted form
\begin{equation}
P_t=ct +\sum_{i=1}^{N_t}\xi_i,
\end{equation}
where $\xi_i$ form an i.i.d. sequence of random variables independent of the counting Poisson process $N$ with intensity $\lambda$; in such a representation $q(dt,dx)=\lambda dt F(dx)$, where $F$ is the distribution of $\xi_i$. Of course, in the case considered here such a representation is not possible in general. The models where $\Pi_P$ is not a finite measure are interpreted as a convenient asymptotic approximation in the case of a large number of small claims. They were considered in the literature starting from the seminal paper by Paulsen \cite{Paul-93}, where a much more general model of the generalized Ornstein--Uhlenbeck process (GOU) was suggested; more recent results and further references can be found in \cite{KP2020}.

We consider on $\bbr_+$ the process $X=X^{u}$ with
\begin{equation}
X^{u}_t=u+\int_{]0,t]}X^{u}_{s}\,dR_s+P_t,
\end{equation}
or, in the differential form,
\begin{equation}
dX^{u}_t=X^{u}_{t} (adt+\sigma d W_t)+ c dt + \int x\mu (dt,dx), \qquad X^{u}_0=u.
\end{equation}

In the actuarial context, $X^{u}$ describes the evolution of the capital of an insurance company starting its activity at time zero with the initial capital $u>0$ and investing its entire reserve in a risky asset whose price evolves as a geometric Brownian motion.

Let $\tau^u:=\inf \{t\colon X^u_t\le 0\}$ be the exit time of the process $X^u$ from $(0,\infty)$. In the actuarial context, it is called the instant of ruin, and the functions $\Psi (u):=\P[\tau^u<\infty]$ and $\Phi(u):=1-\Psi (u)$ are referred to as the ruin and survival probabilities, respectively. It is natural to extend these definitions to the whole real line by putting $\Psi(u)=1$ and $\Phi(u)=0$ for $u\le 0$. Note that $\Psi$ and $\Phi$ are monotone (hence, continuous except on a countable set and differentiable almost everywhere) and may have a discontinuity at zero.

It is plausible to guess that with such extensions both functions will satisfy the stationary Fokker--Planck-type equation $\cL G(u)=0$, which in the case considered looks like
\begin{equation}
\label{eqG}
\frac 12 \sigma^2u^2 G''(u)+ (au+c) G' (u)+\int ( G(u+z)- G(u))\Pi_P(dz)=0.
\end{equation}
Indeed, it is easy to check that they satisfy this equation in the viscosity sense. An appropriate boundary value problem can be investigated, see \cite{Bel-Kab}. For a more advanced framework of the generalized Ornstein--Uhlenbeck process, the Fokker--Planck equation was obtained in \cite{AK2025}. At the moment, it is not clear how to use viscosity solutions for asymptotic analysis. Another approach consists in proving a priori smoothness of the ruin probabilities and then get the equation (\ref{eqG}). Until recently, the arguments were purely probabilistic and available only for the case where $P$ is a compound Poisson process with drift, that is, when $\Pi_P$ is finite. These arguments require restrictive conditions, like the existence of a smooth density of the jump distribution. The results were obtained using rather delicate and cumbersome estimates, see \cite{KP2016}.

The method developed in \cite{Prom} and \cite{AK2026} reverses the order of arguments. The idea can be traced back to the works \cite{Grandits2004} and \cite{Bel2014}. One can use the arguments in the spirit of the ``verification lemma'' of stochastic optimal control theory. If there is a smooth solution of the Bellman equation with some additional properties, then this solution can be identified with the Bellman function. In the present context one needs only to find a positive solution of (\ref{eqG}) smooth on $(0,\infty)$ and such that $G(\infty)=1$. A straightforward use of the It\^o formula allows one to conclude that it is the survival probability function.

In this paper {\bf we assume} that the measure $\Pi_P$ has no atoms, $\Pi_P(\bbr)=\infty$, and there is $p\in (0,1)$ such that $\Pi_P(|h|^{p})<\infty$. It follows that $\Pi_P(|h|^{p'})<\infty$ for all $p'\ge p$. Since $p$ can be arbitrarily close to one, this means that we need only slightly stronger integrability than that required in the model formulation, i.e., $\Pi_P(|h|)<\infty$. However, in the model of annuity payments, we require $p\in (0,1/2)$.

{\sl Notations.} We use, where it is more convenient, either the parameter $\beta:=2a/\sigma^2 - 1$ or $\gamma:=\beta+1=2a/\sigma^2$ (both frequently used in the existing literature) and {\bf we assume} that $\beta>0$ (that is, $\gamma>1$). The parameter $\beta $ makes formulae more readable in expressions involving the ruin probability, while $\gamma$ suits better for the derivative of the latter.

Recall that the case where $\Pi_P(\bbr) < \infty$, already studied in \cite{Prom} and \cite{AK2026}, corresponds to the model where $P$ is a compound Poisson process.

We start with a verification result. It asserts that the positive solution of the equation (\ref{eqG}) tending to one at infinity coincides with the survival probability.

\begin{proposition}
\label{prop}
Suppose that an increasing function $G: \bbr\to [0,1]$ is such that $ G(u)=0$ for $u\leq0$, $G$ is twice continuously differentiable on $(0,\infty)$ and the equation (\ref{eqG}) holds. If $G(\infty)=1$, then $G=\Phi$.
\end{proposition}
\begin{proof}
Let us apply It\^o's formula to $G(X^u_{t\wedge \tau^u})$ where $G$ is the classical solution of (\ref{eqG}) with $G(\infty)=1$. As a result, we get that $G(X^u_{t\wedge \tau^u})-G(u)=M_t$, where
\[
M_t=\int_0^{t\wedge \tau^u}G'(X^u_{s\wedge \tau^u}) \sigma X^u_s dW_s+\int_0^{t\wedge \tau^u}\int (G(X^u_{s-}+x)-G(X^u_{s-}))(\mu-q)(ds,dx).
\]
The process $(M_t)$ is a bounded martingale starting from zero. Hence, $\E [M_\tau]=0$ for any stopping time $\tau$.

Note that the assumption $\gamma >1$ ensures that $X^u_{t}\to \infty$ as $t\to \infty$ (a.s.) on the set $\{\tau^u=\infty\}$, that is, $X^u_{\infty}=\infty$ (a.s.) on this set. Indeed, by the stochastic version of the Cauchy formula $X^u_t=S_tY_t^u$, where
\[
S_t:=\exp\{( a- \sigma^2/2)t + \sigma W_t\} =\exp\{(\sigma^2/2)t( \beta +2W_t/(\sigma t))\}, \quad Y_t^u:=u+\int_0^tS_s^{-1}\,dP_s.
\]
By the law of large numbers, $W_t/t\to 0$ (a.s.) as $t\to \infty$ and, therefore, $S$ is exponentially increasing. On the set $\{\tau^u=\infty\}$ the process $Y$ is strictly positive. Indeed,
\[
Y_t^u = u + \int_{0}^{t} S_s^{-1} dP_s \to Y_{\infty}^u := u + c \int_{0}^{\infty} S_s^{-1} ds + \int_{0}^{\infty} S_s^{-1} \int x \mu(ds, dx), \quad t \to \infty.
\]
Put $\varepsilon > 0$. Let $\xi^{\varepsilon}_1$ be the first jump larger than $\varepsilon$ and $T^{\varepsilon}_1$ be the time of the first jump. Due to the fact that the measure $\Pi_{P}$ has no atoms, the random variable $\xi^{\varepsilon}_1$ has a continuous distribution. Let $\zeta^{\varepsilon} := S_{T^{\varepsilon}_1} (Y_{\infty}^u - S^{-1}_{T^{\varepsilon}_1} \xi^{\varepsilon}_1)$. Due to the independence of $W$ and $P$, we get $\xi^{\varepsilon}_1 \neq \zeta^{\varepsilon}$ (a.s.) and $Y_{\infty}^u > 0$ on $\{ \tau^u = \infty\}$. Hence, on this set $X^u_{\infty}=\infty$. It follows that
\[
G(u)= \E[ G(X^u_{\tau^u})I_{\{\tau^u=\infty\}}]+\E[ G(X^u_{\tau^u})I_{\{\tau^u<\infty\}}]=G(\infty)\P(\tau^u=\infty)=\Phi(u)
\]
and we get the result.
\end{proof}

Our assumption about the function $G$ implies that $G(u)=G(0+)+\int_0^u g(s)\,ds$ for $u>0$, where the density $g$ is continuously differentiable on $(0,\infty)$, and that $g$ satisfies on $(0,\infty)$ the equation
\begin{equation}
\label{eq-g}
\frac 12 \sigma^2u^2 g'(u)+ (au+c) g (u)+\int ( G(u+z) - G(u))\,\Pi_P(dz)=0.
\end{equation}

We treat this equation separately for two models:

{\bf Model 1 (Non-life insurance).} The constant $c>0$ (upward trend), $\Pi_P(\bbr_+)=0$ (downward jumps).

{\bf Model 2 (Annuity payments).} The constant $c<0$ (downward trend), $\Pi_P(\bbr_-)=0$ (upward jumps).

Our main results look similar for both models, but the conditions are different.

\begin{theorem}[Non-life insurance]
\label{main1}
Suppose that $\Pi_P(|\bar h|^{\beta})<\infty$ and there is $p\in (0,1)$ such that $\Pi_P(|h|^{p})<\infty$. Then $u^{\beta}\Psi(u)\to C$ as $u\to \infty$, where $C>0$ is a constant.
\end{theorem}

\begin{theorem}[Annuity payments]
\label{main2}
Suppose that $\Pi_P(|\bar h|)<\infty$ and there is $p\in (0,1/2)$ such that $\Pi_P(|h|^p)<\infty$. Then $u^{\beta}\Psi(u)\to C$ as $u\to \infty$, where $C>0$ is a constant.
\end{theorem}

Integration by parts shows that the first hypothesis of Theorem \ref{main2} admits the equivalent form
\begin{equation}
\label{int1}
\int_1^\infty \Pi_P([u,\infty))\,du <\infty.
\end{equation}
Since $|h|\le |h|^p$ for $p\in [0,1]$, and $\Pi_P(|h|^p)<\infty$, this first hypothesis can be replaced by the following one:
\begin{equation}
\label{int0}
\int_0^\infty \Pi_P([u,\infty))\,du <\infty.
\end{equation}

If $P$ is a compound Poisson process with drift (i.e., the measure $\Pi_P$ is finite), the situation is much simpler. The first hypotheses of the above theorems mean, respectively, that $\E[|\xi_i|^\beta]<\infty$ and $\E[|\xi_i|]<\infty$, that is, they are the same as in \cite{AK2026} and \cite{Prom}. The second hypotheses are fulfilled automatically in such a case.

\section{Model 1}

\begin{theorem}
\label{M1}
Let $\Pi_P(|\bar h|^\beta)<\infty$ where $\beta>0$. Suppose that there is $p\in (0,1)$ such that $\Pi_P(|h|^{p})<\infty$. Then there exists an increasing function $G: \bbr\to [0,1]$ with $ G(u)=0$ for $u\le 0$, the restriction $G\vert _{(0,\infty)}\in C^2((0,\infty))$, and (\ref{eqG}) holds on $(0,\infty)$.
\end{theorem}
\begin{proof}
To conform to the notation used in the literature, we put $\alpha:=2c/\sigma^2$. We also introduce the measure $F(dz)$ by putting $F(\Gamma):=(2/\sigma^2)\Pi_P(-\Gamma)$, where $\Gamma$ is a Borel set in $\bbr$. With this notation, the equation (\ref{eqG}) can be written as
\begin{equation}
\label{red1}
u^2 G''(u)+ (\gamma u+\alpha) G' (u)+\int_{0}^{\infty} ( G(u-z)- G(u))\,F(dz)=0.
\end{equation}
The equation (\ref{eq-g}) can be transformed into the following one:
\begin{equation}
\label{red2}
u^2 g'(u)+ (\gamma u+\alpha) g(u) =\xi \bar F_u+\bar F_u\int_0^u g(s)\,ds+\int_0^u \int_{u-z}^u g(s)\,ds\,F(dz),
\end{equation}
where the constant $\xi:=G(0+)$ and $\bar F_u:=F([u,\infty))$. Note that $\bar F_u\uparrow \infty$ as $u\downarrow 0$.

Let us denote the second and the third terms on the right-hand side of (\ref{red2}) by $A_1g(u)$ and $A_2 g(u)$, that is, we put
\[
A_1g(u):=\bar F_u\int_0^u g(s)\,ds, \qquad A_2 g(u):=\int_0^u \int_{u-z}^u g(s)\,ds\,F(dz).
\]
The left-hand side of the equation (\ref{red2}) admits an integrating factor $u^\gamma e^{-\alpha/u}$, which allows one to rewrite (\ref{red2}) in the form
\begin{equation}
\label{red3}
(u^\gamma e^{-\alpha/u}g(u))'=u^{\gamma-2}e^{-\alpha/u}(\xi \bar F_u+A_1g(u)+A_2g(u)).
\end{equation}
Necessarily, $u^\gamma e^{-\alpha/u}g(u)\to 0$ as $u\downarrow 0$, and the above integro-differential equation can be represented as the integral equation
\[
g(u)= \xi u^{-\gamma} e^{\alpha/u}\int_0 ^u t^{\gamma-2}e^{-\alpha/t}\bar F_t\,dt+ T g(u),
\]
where $T$ is the linear mapping defined by
\[
T g(u):=u^{-\gamma} e^{\alpha/u}\int_0^u t^{\gamma-2}e^{-\alpha/t}(A_1g(t)+A_2g(t))\,dt.
\]

We expect that the solution of this integral equation has a singularity at zero. To treat it, we introduce the ``normalized'' function $g^o (u)=g^ o(u,p)=u^{p}g(u)$ and reformulate the equation in terms of $g^o$. We get that
\begin{equation}
\label{norm}
g^o(u)= \xi u^{p-\gamma} e^{\alpha/u}\int_0 ^u t^{\gamma-2}e^{-\alpha/t}\bar F_t\,dt+ u^{p}T(u^{-p} g^o(u)),
\end{equation}
or, in the operator form, $g^ o(u)=q(u)+\cT g^o(u)$, where $q(u)$ is the known function given by the first term on the right-hand side and $\cT g^o(u)$ is the abbreviation for the second term. We show that the integrability of $|h|^p$ implies that this equation admits, on a sufficiently small interval $[0,u_0]$, a unique continuous solution for any $\gamma>1$. We formulate this claim as follows.

\begin{lemma}
Suppose that $ F( |h|^p)<\infty$ for some $p\in (0, 1)$. Then there is $u_0\in (0,1)$ such that the integral equation (\ref{norm}) admits a unique solution $g^ o(u,p)\in C([0,u_0])$.
\end{lemma}
\begin{proof}
Note that $\sup_{t\le 1} t^p\bar F_t=:C_F$ is finite. Indeed, for any $t\le 1$
\begin{equation}
\label{bound}
\infty>M_p:=\int_0^1 z^p F(dz)\ge \int_t^1 z^p F(dz)\ge t^p F([t,1))\ge t^p\bar F_t-\bar F_1.
\end{equation}
Recall that
\begin{equation}
\label{explicit}
e^{\alpha/u}\int_0^u t^{-2}e^{-\alpha/t}\,dt=\frac{1}{\alpha}.
\end{equation}
Taking into account that $t^{\gamma-p}\le u^{\gamma-p}$ for $t\in [0,u]$, we get that for $u\le 1$
\[
u^{p-\gamma} e^{\alpha/u}\int_0^u t^{\gamma-2}e^{-\alpha/t}\bar F_t\,dt=u^{p-\gamma} e^{\alpha/u}\int_0^u t^{\gamma-p-2}e^{-\alpha/t}t^p\bar F_t\,dt\le C_F/\alpha <\infty.
\]
Applying, for $t\le u$, the bounds
\[
\int_0^ts^{-p}g^o(s)\,ds\le \frac{1}{1-p}t^{1-p}||g^o||_u, \qquad t\bar F_t\le u^{1-p}C_F,
\]
where $||g^o||_u:=\sup_{s\le u}|g^o(s)|$, and combining them with the above arguments, we get that
\[
u^{p-\gamma} e^{\alpha/u}\int_0^u t^{\gamma-2}e^{-\alpha/t}\bar F_t\int_0^ts^{-p}g^o(s)\,ds\,dt \le\frac{1}{\alpha(1-p)}C_F u^{1-p}||g^o||_u.
\]

For $x\in [0,1]$ and $a\in (0,1)$ we have the inequality $1 - (1-x)^a\le x$. Using it, we get:
\begin{eqnarray*}
&&u^{p - \gamma} e^{\alpha/u}\int_0^u t^{\gamma-2}e^{-\alpha/t}\int_0^t (t^{1-p}-(t-z)^{1-p})\, F(dz)\,dt\\
& = & u^{p - \gamma} e^{\alpha/u}\int_0^u t^{\gamma-1-p}e^{-\alpha/t}\int_0^t (1-(1-z/t)^{1-p})\, F(dz)\,dt\\
&\le & F(hI_{[0,u]}) \left[ u^{p - \gamma} e^{\alpha/u}\int_0^u t^{\gamma - p - 2}e^{-\alpha/t}\,dt \right] \le F(hI_{[0,u]})/\alpha.
\end{eqnarray*}

It follows that
\[
u^{p-\gamma} e^{\alpha/u}\int_0^u t^{\gamma-2}e^{-\alpha/t}\int_0^t \int_{t-z}^t s^{-p} g^o(s)\,ds\, F(dz)\,dt \le \frac{1}{\alpha(1-p)}F(hI_{[0,u]}) ||g^o||_u.
\]
The integrability of $h$ implies that $F(hI_{[0,u]}):=\int_0^uzF(dz)\downarrow 0$ as $u\downarrow 0$.

The above estimates lead to the conclusion that the norm of $\cT$, considered as a continuous linear operator on the space $C([0,u_0])$, is less than one when $u_0$ is sufficiently small. That is, the mapping $\cT$ is a contraction for sufficiently small $u_0$ and so is the affine mapping $f\mapsto q+\cT f$. We get the result by applying the Banach fixed point theorem.
\end{proof}

Using the fixed point $g^ o$ provided by the lemma, we put $g_*(u):=u^{-p}g^o (u)$. It is a positive function continuous on $(0,u_0]$, having an integrable singularity at $0$, and solving the equation (\ref{red3}) on this interval. In other words, we established the local existence and uniqueness for this equation near the origin. Now we show that appropriate additional assumptions allow us to extend $g_*(u)$ to $(0,\infty)$ and control the behavior of such an extension at infinity.

To treat the convergence to zero of solutions at infinity, it is helpful to use the normalization trick, but this time with an amplifying factor. To this end, we introduce the function $g_{\vartheta}(u):=u^{\vartheta } g(u)$, where $\vartheta> 1$.

Integrating (\ref{red3}) from $u_0$ upward, we obtain the representation
\begin{equation}
\label{gg}
g(u)u^{\gamma}e^{-\alpha/u}= Y(u)+ \int_{u_0}^u t^{\gamma-2}e^{-\alpha/t} \int_{u_0}^t g(s) \bar F_{t-s}\,ds\,dt
\end{equation}
where we isolated the already known part $g_*$ of $g$ denoted by
\begin{eqnarray*}
Y(u)&:=&u_0^\gamma e^{-\alpha/u_0} g_*(u_0)+ \int_{u_0}^u t^{\gamma-2}e^{-\alpha/t} v(t)\,dt
\end{eqnarray*}
with the abbreviation
\[
v(t):=\xi \bar F_t+\bar F_t\int_0^{u_0} g_*(s)\,ds+\int_0^{u_0} g_*(s) F([t-s,t])\,ds=\xi \bar F_t+\int_0^{u_0} g_*(s) \bar F_{t-s}\,ds.
\]

Introducing the kernel
\begin{equation}
\label{Kus}
K_\vartheta(u,s):=s^{- \vartheta } \int_{s}^u t^{\gamma-2}e^{-\alpha/t} \bar F_{t-s}\,dt,
\end{equation}
we rewrite (\ref{gg}) as the integral equation for $g_{\vartheta}(u)$:
\begin{equation}
\label{eq-tilde}
g_{\vartheta}(u)= u^{\vartheta - \gamma } e^{\alpha/u}Y(u)+ u^{\vartheta - \gamma } e^{\alpha/u}\int_{u_0}^u g_{\vartheta}(s) K_\vartheta(u,s)\,ds.
\end{equation}
Note that
\[
\bar F_{t-s}=\bar F_{t-s}I_{\{t-s\le 1\}}+ \bar F_{t-s}I_{\{t-s>1\}}\le C_F(t-s)^{-p}I_{\{t-s\le 1\}}+\bar F_1.
\]
It follows that the integral in the definition of $K_\vartheta(u,s)$ is finite for all $s>u_0$ and that the kernel $K_\vartheta(u,s)$ is continuous. Thus (\ref{eq-tilde}) is a Volterra equation of the second kind. It admits a continuous solution $g^{*}_{\vartheta}$ on $[u_0,\infty)$ which is unique, see \cite{Burton}, Th. 2.1.1. Since \\ $C_F:=\sup_{t\le 1} t^p\bar F_t<\infty$, we easily derive from (\ref{gg}) that $g$ has a continuous derivative.

Note that there is a freedom in the choice of the point $u_0$ ensuring that our construction leads to a unique solution of (\ref{red3}) continuously differentiable on $(0,\infty)$.

\begin{lemma}
Let $\vartheta \in (1, \gamma ]$. Suppose that
\begin{equation}
\label{J}
J_\vartheta (u):= \int_{u}^\infty t^{\vartheta - 2}\bar F_t\,dt <\infty
\end{equation}
for some $u>0$ (hence, for all $u>0$). Then there exists a constant $C_{\vartheta}$ such that $g_{\vartheta}(u) \le C_{\vartheta}$ on $[u_0,\infty)$. Furthermore, if $\vartheta = \gamma $, then there is a constant $C_\infty >0$ such that $ g_{\vartheta}(u) \to C_\infty $ as $u \to \infty$.
\end{lemma}

\begin{proof}
Without loss of generality, we assume that $u_0<1$. First, we check that the function $u^{\vartheta - \gamma } e^{\alpha/u}Y(u)$ is bounded on $[u_0,\infty)$.

The assumption $\vartheta\le \gamma$ implies the inequality $u^{\vartheta-\gamma } e^{\alpha/u}\le t^{\vartheta-\gamma } e^{\alpha/t}$ for $t\le u$. Using this observation, we get that
\[
u^{\vartheta-\gamma } e^{\alpha/u} \int_{u_0}^u t^{\gamma-2} e^{-\alpha/t}\bar F_t\,dt \le \int_{u_0}^{u} t^{\vartheta - 2} \bar F_t \,dt \le J_\vartheta (u_0)<\infty.
\]
On the interval $[0,u_0]$ we have the representation $g_*(u)=u^{-p}g^o(u)$ where the function $g^o(u)$ is continuous, hence, bounded. Recalling that $C_F:=\sup_{v\le 1}v^p\bar F_v<\infty $, we get that
\begin{eqnarray*}
&&\int_{u_0}^u t^{\vartheta-2} \int_0^{u_0}g_*(s) \bar F_{t-s} I_{\{t-s\le 1\}}\,ds\,dt\\
&&\qquad\le C_F ||g^o||_{u_0} \int_{u_0}^u t^{\vartheta-2} \int_0^{u_0}(t-s)^{-p} s^{-p} I_{\{t-s\le 1\}} \,ds\,dt
\end{eqnarray*}
and the integral on the right-hand side is dominated by a constant:
\begin{eqnarray*}
&&\int_{u_0}^u t^{\vartheta-2} \int_0^{u_0}(t-s)^{-p} s^{-p} I_{\{t\le u_0+1\}} \,ds\,dt\\
&&\qquad\le \max\{u_0^{\vartheta-2},(u_0+1)^{\vartheta-2}\} \int_0^{u_0}(u_0-s)^{-p} s^{-p}\,ds<\infty.
\end{eqnarray*}

By Fubini's theorem,
\[
\int_{u_0}^u t^{\vartheta-2} \int_0^{u_0} g_*(s) \bar F_{t-s} I_{\{t-s>1\}} \,ds\,dt = \int_{0}^{u_0} g_*(s) \int_{u_0}^{u} t^{\vartheta - 2} \bar F_{t-s}I_{\{t-s> 1\}}\, dt\, ds.
\]
To check that the right-hand side of this equality is also bounded, it is sufficient to verify that for $s\le u_0$ the integral in $dt$ is bounded uniformly in $s$. We have:
\[
\int_{u_0}^{\infty} t^{\vartheta - 2} \bar F_{t-s}I_{\{t-s>1\}}\, dt \le \int_{1}^{\infty} t^{\vartheta - 2}\bar F_{t-u_0}dt<\infty,
\]
because the inequality $t-s>1$ implies that $t>1$, and
\begin{eqnarray*}
\int_{1}^{\infty} t^{\vartheta - 2}\bar F_{t-u_0}dt &=&\int_{1-u_0}^{\infty} (1+u_0/r)^{\vartheta - 2}r^{\vartheta-2} \bar F_{r}\,dr\\
&\le& \max\{1,(1-u_0)^{2-\vartheta}\}\, J_\vartheta (1-u_0)<\infty
\end{eqnarray*}
by virtue of the integrability assumption (\ref{J}) and our convention that $u_0<1$.

It follows from the above estimates that on the interval $[u_0,\infty)$ the continuous positive function $u^{ \vartheta -\gamma } e^{\alpha/u}Y(u)$ is bounded by some constant $C_Y>0$.

Arguing in the same spirit, we deduce from (\ref{Kus}) that
\begin{eqnarray*}
u^{\vartheta-\gamma } e^{\alpha/u} K_\vartheta(u,s)&\le & e^{\alpha/u_0}V_{\vartheta }(s), \qquad V_{\vartheta }(s):=s^{-\vartheta } \int_{s}^{\infty} t^{\vartheta - 2}e^{-\alpha/t} \bar F_{t-s}\,dt.
\end{eqnarray*}
We represent  the integral as a sum of two integrals and estimate them separately.   For $s\ge u_0$ we have:  
\begin{eqnarray*}
\int_{s}^{\infty} t^{\vartheta - 2} \bar F_{t-s}I_{\{t-s\le 1\}}\,dt&\le& C_F \int_{s}^{s+1} t^{\vartheta - 2} (t-s)^{-p}\,dt\\
&=& C_F s^{\vartheta - p-1}\int_{s/(s+1)}^{1} r^{p-\vartheta} ( 1-r)^{-p}\,dr\\
&\le & C_F s^{\vartheta - p-1}\int_{u_0/(u_0+1)}^{1} r^{p-\vartheta} ( 1-r)^{-p}\,dr.
\end{eqnarray*}
The change of variable $r=t-s$ gives: 
\[
\int_{s}^{\infty} t^{\vartheta - 2} \bar F_{t-s}I_{\{t-s> 1\}}\,dt=\int_{1}^{\infty} (r+s)^{\vartheta - 2} \bar F_{r}\,dr.
\]
If $\vartheta\le 2$, then $(r+s)^{\vartheta-2}\le r^{\vartheta-2}$ and the right-hand side of the above identity is dominated by $J_\vartheta(1)$. If $\vartheta>2$, then $(r+s)^{\vartheta-2}\le 2^{\vartheta-2}(r^{\vartheta-2}+s^{\vartheta-2})$ and the right-hand side does not exceed $2^{\vartheta-2}(J_\vartheta(1)+s^{\vartheta-2}J_2(1))$, where $J_2(1)\le J_\vartheta(1)$.
It follows that there is a constant $C_V>0$ such that $V_\vartheta(s)\le C_V (s^{-p-1}+s^{-\min(\vartheta,2)})$, and, therefore, it is  integrable on $[u_0,\infty)$ because $\vartheta>1$.

Using the above estimates we get, for some constant $C$, that for $u\ge u_0$
\[
g_{\vartheta}(u) \leq C + e^{\alpha/u_0}\int_{u_0}^{u} g_{\vartheta}(s)V_{\vartheta }(s)\, ds
\]
and by the Gronwall--Bellman lemma
\[
g_{\vartheta}(u)\le C e^{e^{\alpha/u_0}\int_{u_0}^u V_\vartheta(s)\,ds}\le C e^{e^{\alpha/u_0}\int_{u_0}^\infty V_\vartheta(s)\,ds}.
\]
Thus, $C_{\vartheta}:=\sup_{u\ge u_0}g_{\vartheta}(u)<\infty$.

When $\vartheta=\gamma$, the equation (\ref{eq-tilde}) reduces to the following one:
\[
g_{\gamma}(u)= e^{\alpha/u}Y(u)+ e^{\alpha/u}\int_{u_0}^u g_{\gamma}(s) K_\gamma(u,s)\,ds
\]
where $K_\gamma(u,s)\to V_\gamma(s)$ as $u\to \infty$. Therefore,
\[
g_{\gamma}(\infty)=Y(\infty)+\int_{u_0}^{\infty} g_{\gamma}(s)V_\gamma(s)\,ds=C_\infty,
\]
and the lemma is proven.
\end{proof}

Theorem \ref{M1} follows from the above lemma.
\end{proof}

\noindent
\begin{remark} 
{\rm
In the case where $P$ is a compound Poisson process and $F$ is the common distribution of the absolute values of jumps $\xi_i$, the condition $J_\gamma(u) <\infty$ means that $\E |\xi_i|^\beta<\infty$.

Note that the obtained solution depends on the parameter $\xi=G(0+)$. To ensure that the function $G(u)=G(u;\xi)$ is a candidate for the survival probability $\Phi(u )$, we must determine this parameter from the equation
\[
\xi+\int_0^\infty g(s;\xi)\,ds =1.
\]
With such $\xi$, which always exists due to linearity, we obtain Theorem \ref{main1} from Proposition \ref{prop} and Theorem \ref{M1}. }
\end{remark}

\section{Model 2}

In the case of upward jumps, we put $\alpha:=-2c/\sigma^2$, $F(du):=(2/\sigma^2)\Pi_P(du)$. Now $\alpha>0$ and $F(\bbr_-)=0$. The equation (\ref{eqG}) can be written as
\begin{equation}
\label{red1up}
u^2 G''(u)+ (\gamma u-\alpha) G' (u)+\int_0^\infty ( G(z+u)- G(u))\,F(dz)=0
\end{equation}
and the equation (\ref{eq-g}) as
\begin{equation}
\label{red2up}
u^2 g'(u)+ (\gamma u-\alpha) g (u)+\int_0^\infty \int_{u}^{z+u}g(s)\,ds\,F(dz)=0.
\end{equation}
Applying Fubini's theorem, we obtain that the double integral above is equal to
\[
Ag(u):= \int_u^\infty g(s)\bar F_{s-u}\,ds.
\]
Thus, the equation (\ref{red2up}) can be written as
\begin{equation}
\label{red3up}
u^2 g'(u)+ (\gamma u-\alpha) g (u)=-Ag(u),
\end{equation}
or, equivalently,
\begin{equation}
\label{ide}
(u^\gamma e^{\alpha/u}g(u))'=-u^{\gamma-2} e^{\alpha/u}Ag(u).
\end{equation}

\begin{lemma}
Suppose that
\begin{equation}
\label{intfin}
J_2(0):=\int_0^\infty \bar F_s\,ds<\infty.
\end{equation}
Then there exists $u_0\in [1,\infty)$ such that the equation (\ref{ide}) has on $[u_0,\infty)$ a unique solution $g(u)\ge 0$ such that $u^\gamma g(u)\to 1$ as $u\to \infty$.
\end{lemma}
\begin{proof}
Put $\tilde g(u):=u^\gamma e^{\alpha/u}g(u)$. Integrating from $u\ge u_0$ to $\infty$, we obtain from (\ref{ide}) the integral equation for the function $\tilde g(u)$ with $\tilde g(\infty)=1$:
\begin{equation}
\label{ie}
\tilde g(u)=1+\int_u^\infty t^{\gamma-2} e^{\alpha/t} \int_t^\infty s^{-\gamma} e^{-\alpha/s} \tilde g(s)\bar F_{s-t}\,ds\,dt.
\end{equation}
Let us denote by $T\tilde g(u)$ the integral on the right-hand side above and put $||v||^{u_0}:=\sup_{u\ge u_0}|v(u)|$. Then for $u\ge u_0$
\begin{eqnarray*}
T\tilde g(u)&\le& e^{\alpha/u_0}||\tilde g||^{u_0}\int_u^\infty t^{\gamma-2} \int_t^\infty s^{-\gamma} \bar F_{s-t}\,ds\,dt \\
&\le& ||\tilde g||^{u_0} \frac{e^{\alpha/u_0}}{u_0}\int_0^\infty \bar F_s\,ds .
\end{eqnarray*}
Since $s\mapsto \bar F_s$ is an integrable function, $||T\tilde g||^{u_0}\le \theta(u_0)||\tilde g||^{u_0}$, where $\theta(u_0)\downarrow 0$ as $u_0 \uparrow \infty$. Thus, for sufficiently large $u_0$, the linear operator $T$ in the space $C_b([u_0,\infty))$ of bounded continuous functions is a contraction. It follows that the mapping $\tilde g\mapsto 1+T\tilde g$ is also a contraction in this space. By the Banach fixed point theorem it has a fixed point $\tilde g(u)$ which is unique.
\end{proof}

At the second step we analyze the equation (\ref{ide}) considering that on $[u_0,\infty)$ the function $g$ is known and coincides with $g^*(u):= u^{-\gamma}e^{-\alpha/u}\tilde g(u;1)$. Integrating now from $u$ to $u_0$ and multiplying both sides of the obtained equality by $ u^{-\gamma} e^{-\alpha/u}$, we get that
\begin{eqnarray*}
g(u)&=& u^{-\gamma} e^{-\alpha/u} u_0^\gamma e^{\alpha/u_0}g^*(u_0)+ u^{-\gamma} e^{-\alpha/u} \int_{u}^{u_0} t^{\gamma-2}e^{\alpha/t}H(t)\,dt\\
&&+ u^{-\gamma} e^{-\alpha/u}\int_u^{u_0}t^{\gamma-2}e^{\alpha/t}\int_t^{u_0}g(s)\bar F_{s-t}\,ds\,dt,
\end{eqnarray*}
where
\[
H(t):=\int_{u_0}^\infty g^*(s)\bar F_{s-t}\,ds.
\]
Let us denote by $f(u)$ the sum of the first two terms on the right-hand side of the equation above. Changing the order of integration in the third term, we arrive at the integral equation
\begin{equation}
\label{Volterra}
g(u)=f(u)+\int _u^{u_0} K(u,s)g(s)\,ds, \qquad u\in (0,u_0),
\end{equation}
where
\[
K(u,s):=u^{-\gamma} e^{-\alpha/u}\int _u^s t^{\gamma-2}e^{\alpha/ t}\bar F_{s-t}\,dt.
\]
The equation (\ref{Volterra}) is a Volterra equation of the second kind and has a unique continuous solution on $(0,u_0)$ which we denote by $g_*$.

We show that the normalized function $g^o_*(u):=u^{p}g_*(u)$ is bounded near zero. Such a property will imply that $g_*$ is integrable at zero. 
The equation (\ref{Volterra}) is equivalent to the equation
\begin{equation}
\label{Volterra+}
g^o_*(u)=u^{p}f(u)+u^{p}\int _u^{u_0} s^{-p}K(u,s)g^ o_*(s)\,ds, \qquad u\in (0,u_0).
\end{equation}

First, we observe that by L'H\^opital's rule
\[
\lim_{u\downarrow 0} u^{-\gamma} e^{-\alpha/u} \int_{u}^{u_0} t^{\gamma-2}e^{\alpha/t}H(t)\,dt= -\lim_{u\downarrow 0}\frac{u^{\gamma-2}e^{\alpha/u}H(u)}{(\gamma u-\alpha )u^{\gamma-2}e^{\alpha/u}}=\frac{H(0)}{\alpha}<\infty.
\]
Hence, $f$ is bounded on $[0,u_0]$. On this interval $u^{p}f(u)\le \kappa $ where $\kappa$ is a constant. Also, for $s\ge u>0$,
\begin{eqnarray*}
u^{p}s^{-p}K(u,s)&=&u^{p}s^{-p}u^{-\gamma} e^{-\alpha/u}\int _u^s t^{\gamma-2}e^{\alpha/ t}\bar F_{s-t}\,dt\\
&= &u^{p}s^{-p}u^{-\gamma}   e^{-\alpha/u} \int _u^s t^{\gamma-2} e^{\alpha/t}(s-t)^{-p}(s-t)^p\bar F_{s-t}\,dt\\
&\le& \kappa_1  u^{-\gamma} e^{-\alpha/u}\int _u^s t^{\gamma-2} e^{\alpha/t}(s-t)^{-p}\,dt=:\kappa_1V(u,s),
\end{eqnarray*}
where $\kappa_1:= \sup_{t\le u_0}t^p\bar F_t $.

%

\begin{lemma}
\label{uniform}
Let $p\in (0,1)$. Then there is 
a constant $C_0$ such that $V(u,s)\le C_0s^{-2p}$ for all $u$ and $s$ such that  $0<u\le s\le u_0$. 
\end{lemma}
\begin{proof}
The change of variable $y=1/t$ results in the expression 
$$
V(u,s)= \int _{1/s}^{1/u} \Big(\frac 1{uy}\Big)^{\gamma} e^{\alpha (y-1/u)}(s-1/y)^{-p}\,dy.
$$
Let us introduce the notation $\delta:=1/u-1/s$.

Using the further change of variable $x=1/u-y$ and the identity 
  $$
  (s-1/(1/u-x))^{-p} =  s^{-2p}(s+(\delta-x)^{-1})^p 
  $$
we get the following formula
$$
V(u,s)=s^{-2p}\int _{0}^{\delta} \Big(1+\frac {ux}{1-ux}
\Big)^{\gamma} e^{-\alpha x}(s+(\delta-x)^{-1})^{p}\,dx=:s^{-2p} I(u,s). 
$$
It remains to check that  $I(u,s)$ is bounded by a constant depending only on $u_0$. 
Note that 
$$
 \Big(1+\frac {ux}{1-ux}\Big)^\gamma\le  
  \Big(1+\frac {ux}{1-u\delta}\Big)^\gamma=(1+sx)^\gamma \le (1+u_0x)^\gamma 
$$
and 
$$
\big(s+(\delta-x)^{-1}\big)^{p}\le s^p+(\delta-x)^{-p} \le  u_0^p+(\delta-x)^{-p}. 
$$
It follows that 
$$
I(u,s)\le u_0^p\int_0^\infty(1+u_0x)^\gamma e^{-\alpha x}\,dx+
\int_0^\delta(1+u_0x)^\gamma e^{-\alpha x}(\delta-x)^{-p}\,dx. 
$$
Put   $\varphi(x):=(1+u_0x)^\gamma e^{-\alpha x}$,  $\varphi^*:=\sup_{x\ge 0}\varphi(x)$, 
$$
A_0:=\int_0^\infty \varphi(x)\,dx, \qquad B(u,s):=\int_0^\delta \varphi(x)(\delta-x)^{-p}\,dx.
$$ 
In this notation the above inequality has the form $I(u,s)\le u_0^p A_0+B(u,s)$. 
For the  integral $B(u,s)$ we have
$$
B(u,s) =\int_0^\delta \varphi(\delta-y)y^{-p}I_{\{y\le 1\}}\,dy+\int_0^\delta \varphi(\delta-y)y^{-p}I_{\{y>1\}}\,dy.
$$
The second integral on the right-hand side does not exceed $A_0$, while the first one is dominated 
by 
$$
\varphi^*\int_0^{\delta\wedge 1}y^{-p}\,dy\le \varphi^*/(1-p). 
$$
Summarizing, we get  the assertion of the lemma  with 
$C_0=(u_0^p+1)A_0+\varphi^*/(1-p)$.
\end{proof}

Thus, for all $u\in (0,u_0]$, we get from the equation (\ref{Volterra+}) and the obtained bounds that 
\[
g^o_*(u)\le \kappa +\kappa_1 C_0\int_u^{u_0} g^o_*(s)s^{-2p}\,ds.
\]
Recall that in the formulation of Theorem \ref{main2} we assume that $p\in (0,1/2)$, which ensures that the function $s^{-2p}$ is integrable at zero.
By the Gronwall--Bellman lemma, $g^ o_*$ is bounded. Therefore, $g_*$ is integrable on $[0,u_0]$. The function $g^*$ is integrable on $[u_0,\infty)$ and behaves at infinity as $u^{-\gamma}$. It follows that $g_*I_{(0,u_0]}+g^*I_{(u_0,\infty)}$ is a function on $\bbr_+$ which is integrable and solves the equation (\ref{ide}). Due to the homogeneity of this equation, we can multiply $g$ by a positive constant $C_{\infty}$ to ensure that the integral of $\hat g:=C_\infty (g_*I_{(0,u_0]}+g^*I_{(u_0,\infty)})$ over $\bbr_+$ is equal to one.

Put $G(u):=\int_0^u \hat g(s)\,ds$. By virtue of Proposition \ref{prop}, the function $G$ coincides with the survival probability $\Phi$.
The proof of Theorem \ref{main2} is complete. 

\begin{remark} {\rm The estimate of Lemma \ref{uniform} is sharp in $p$ in the following sense. There is 
a constant $c_0>0$ such that for any $s\in (0,u_0]$ there exists $u_s\in (0,s)$ with 
$V(u_s,s)\ge c_0 s^{-2p}$. Indeed, take $u_s:=\alpha s/(\alpha+s)$. In such a case $\delta=1/\alpha$ and for $x\in (0,1/\alpha)$ we have the inequalities 
$$
1/(1-u_sx)>1, \qquad e^{-\alpha x}\ge e^{-1}, \qquad (s+(1/\alpha-x)^{-1})^{p}\ge (1/\alpha-x)^{-p}
$$  
implying that 
$$
I(u_s,s)\ge e^{-1}\int_0^{1/\alpha}(1/\alpha-x)^{-p}dx=e^{-1}\alpha^{p-1}/(1-p). 
$$
}
\end{remark}

\section*{Declarations}

\noindent{\bf Author contributions} Both authors contributed equally to the conceptualization, methodology, formal analysis, investigation, and writing of the manuscript. Both authors have read and agreed to the final version of the manuscript.

\medskip
\noindent{\bf Funding} This research received no external funding.

\medskip
\noindent{\bf Data availability} Not applicable.

\medskip
\noindent{\bf Conflict of interest} The authors declare no conflict of interest.

\end{document}